\documentclass[12pt]{article}

\usepackage[left=3cm,right=2cm]{geometry}

\usepackage[utf8]{inputenc}

\usepackage[noadjust]{cite}

\usepackage{makecell}
\usepackage{multirow}
\usepackage{booktabs}
\usepackage{amssymb}
\usepackage{amsmath}
\usepackage{mathtools}
\usepackage{color}
\usepackage{xspace}
\usepackage{xparse}
\usepackage{subfigure}
\usepackage{bbm}
\usepackage{enumerate}
\usepackage{graphicx}

\usepackage{algorithm}
\usepackage{algorithmicx}
\usepackage{algpseudocode}

\usepackage{bm}
\usepackage{tikz}
\usepackage{tkz-euclide}
\usepackage{chngcntr}
\usepackage{verbatim}
\usepackage{mathtools}
\usepackage{esvect}
\usepackage{subfiles}

\usepackage{atbegshi}%
\usepackage{array, tabularx, caption, boldline}
\usepackage{url}
\usepackage[misc,geometry]{ifsym} 
\usepackage{enumitem}
\usepackage[colorinlistoftodos,textsize=scriptsize]{todonotes}

\usepackage[breaklinks=true,colorlinks=true,linkcolor=red, citecolor=blue, urlcolor=blue]{hyperref}%

\def\R{\mathbb{R}}
\newcommand{\E}{{\mathbb E}}

\def\R{\mathbb R}
\def\E{\mathbb E}

\def\e{\varepsilon}
\def\la{\langle}
\def\ra{\rangle}
\def\vp{\varphi}
\def\y{\mathbf{y}}

\def\x{\mathbf{x}}

\def\lm{\lambda}

\def\bld{\boldsymbol}

\def\R{\mathbb{R}}

\def\R{\mathbb R}
\def\E{\mathbb E}

\def\e{\varepsilon}
\def\la{\langle}
\def\ra{\rangle}
\def\vp{\varphi}
\def\y{\mathbf{y}}

\def\x{\mathbf{x}}

\def\lm{\lambda}

\def\bld{\boldsymbol}

\def\pd#1{{\color{blue}#1}} 
\def\dd#1{{#1}} 
\def\dv#1{{#1}}
\def\ga#1{{#1}}
\def\gav#1{{#1}} 
\def\avg#1{{#1}}
\def\sasha#1{{#1}}
\def\g#1{{#1}}

\def\dm#1{{#1}}
\def\dmdv#1{{#1}}

\usepackage[percent]{overpic}

\newtheorem{theorem}{Theorem}

\newtheorem{remark}[theorem]{Remark}

\usepackage[colorinlistoftodos,textsize=scriptsize]{todonotes}

\title{
Decentralized and \g{Parallel} Primal and Dual Accelerated Methods  for Stochastic Convex Programming Problems
}

\author{ Darina Dvinskikh, and Alexander Gasnikov 
    \thanks{D. Dvinskikh
    (\textit{darina.dvinskikh@wias-berlin.de)} is with  Weierstrass Institute for Applied Analysis and Stochastics, and Moscow Institute of Physics and Technology, and Institute for Information Transmission Problems RAS.
    A. Gasnikov (\textit{gasnikov@yandex.ru}) is with Moscow Institute of Physics and Technology, and Institute for Information Transmission Problems RAS, and  Weierstrass Institute for Applied Analysis and Stochastics
        }  
}

\date{}

\begin{document}
\maketitle
\mathtoolsset{showonlyrefs}

\begin{abstract}
We introduce primal and dual stochastic gradient oracle methods for decentralized convex optimization problems. Both {for} primal and dual oracles\dm{,} the proposed methods are optimal in terms of \dd{the number of} communication steps. 
However, \dd{for all classes of the objective, the} optimality in terms of \dd{the number of} oracle calls per node  takes place \dd{only} up to a logarithmic factor and the notion of smoothness. By using mini-batching technique\dm{,} we show that \dm{the} proposed methods with stochastic oracle can be additionally parallelized at each node. \dm{The} \avg{considered algorithms can be applied to many data science problems and inverse problems.}
\end{abstract}

\section{Introduction}

We consider the stochastic convex optimization problem
\begin{equation}
\label{SP}
\min_{x\in Q \subseteq \mathbb{R}^n} f(x) := \E[f(x,\xi)].    
\end{equation}
Such kind of problems arise in many applications of data science 
\cite{shalev2014understanding,shapiro2014lectures} and mathematical statistics \cite{spokoiny2012parametric}. To solve this problem with \gav{the average precision $\e$}   in the function value \dm{(i.e., to find such $x^N$ that $\E f(x^N) - \min\limits_{x\in Q } f(x)\leq \e$)}\dm{,}  one can use  stochastic gradient (mirror) descent \cite{juditsky2012first-order} with
\begin{equation}\label{SNS}
\min\left\{O\left(\frac{M^2R^2}{\e^2}\right),~ O\left(\frac{M^2}{\mu\e}\right)\right\}
\end{equation}
number of calculations of unbiased stochastic subgradients $\nabla f(x,\xi)$. Here   $\mu$  is the constant of strong convexity of $f$,   \dm{ $R = \|x^0 - x_*\|_2$ is the Euclidean distance between starting point $x^0$ and the solution $x^*$ of \eqref{SP} that corresponds to the minimum of this norm}. We also used
$\E[\|\nabla f(x,\xi)\|_2^2]\le M^2$. Generally\dm{,}   we can parallelize \eqref{SNS}  \dm{on no more than} $\tilde{O}(1)$ processors \g{by using batch-parallelization} \cite{dvurechensky2018parallel}. If we additionally assume that $f$ has $L$-Lipschitz (continuous) gradient  and $\E[\|\nabla f(x,\xi) - \nabla f(x) \|_2^2]\le \sigma^2$, then  \eqref{SNS}  \dm{is replaced by}
\begin{equation}\label{SS}
\hspace{-0.4mm}\min\left\{O\left(\sqrt{\frac{LR^2}{\e}}\right)+O\left(\frac{\sigma^2R^2}{\e^2}\right),~ O\left(\sqrt{\frac{L}{\mu}}\ln\left(\frac{\mu R^2}{\e}\right)\right)+O\left(\frac{\sigma^2}{\mu\e}\right)\right\}
\end{equation}
\dv{number of calculations of unbiased stochastic subgradients $\nabla f(x,\xi)$} by using batch parallelization \cite{devolder2013exactness,dvurechensky2016stochastic,gasnikov2018universal,ghadimi2013stochastic}.
{In this case} we can parallelize \dv{subgradients} calculations \ga{\dm{on} no more than}
\begin{align*}
O\left(\frac{{\sigma^2R^2}/{\e^2}}{\sqrt{{LR^2}/{\e}}}\right) ~ \text{ or  } ~ 
O\left(\frac{{\sigma^2}/{(\mu\e)}}{\sqrt{{L}/{\mu}}\ln\left({\mu R^2}/{\e}\right)}\right)
\end{align*}
processors \ga{(depending on where the minimum in \eqref{SS} is reached)}. Notice that this is much better than in previous case. Since this result cannot be improved \cite{woodworth2018graph}, it is the best possible way (in general) to solve \eqref{SP} by using parallel architecture in online context \cite{shalev2009stochastic}.  

\dd{For} many reasons, in some situations in practice\dm{,} it \dd{can} \dv{be} impossible to organize model-based request\footnote{\dv{For desired $x$ and  independently generated $\xi$, a request returns $\nabla f(x,\xi)$. This allows not to keep the set of functions $\{f(x,\xi^k)\}_k$ for different $k$ in the memory.}} for calculating stochastic gradient $\nabla f(x^k,\xi^k)$ in online regime. 
Typically, in machine learning applications \cite{hastie2001elements,shalev2014understanding}\dm{,} instead of online access to $\left\{\nabla f(x^k,\xi^k)\right\}_{k=1}^m$ we have offline access. This means that the set of functions $\left\{f(x,\xi^k)\right\}_{k=1}^m$ are stored  in \dd{the} memory and to use \dd{them} in algorithms\dm{,} we need to request corresponding function and then calculate its gradient. 
This may \dd{significantly} change the complexity of the problem.  Indeed, it is known  from \cite{guigues2017non-asymptotic,shalev2009stochastic,shapiro2014lectures} that with high probability  the exact solution of  problem
\begin{equation}
\label{sum_rand}
 \min_{x\in Q \subseteq \R^n} \tilde{f}(x):=\frac{1}{m}\sum_{k=1}^m f(x,\xi^k)   
\end{equation}
is an $\e$-solution (in the function value) of problem \eqref{SP} if
\[m = \min\left\{\tilde{O}\left(\frac{nM^2R^2}{\e^2}\right),\tilde{O}\left(\frac{M^2}{\mu\e}\right)\right\}.\]
\dv{If $\mu = 0$ or it is small enough one may use a regularization technique (see e.g., \cite{gasnikov2017modern,shalev2009stochastic}). This allows to reduce the first part of the estimate  from $\tilde{O}\left(nM^2R^2/\e^2\right)$ to  $\tilde{O}\left(M^2R^2/\e^2\right).$}
Moreover, we cannot typically find the exact solution of \eqref{sum_rand} but in the $\mu$-strongly convex (or regularized) case it suffices to solve
\gav{\eqref{sum_rand}}
with accuracy  $O(\mu\e^2/M^2)$ (see \cite{shalev2009stochastic}).

To solve \eqref{sum_rand} in offline context we have to store $\{f(x,\xi^k)\}_{k=1}^m$  in the memory. Since $m$ can be large, centralized distributed architecture is often more preferable in this context \cite{bertsekas1989parallel}. In the general case\dm{,} centralized architecture is based on communication network \dv{and it can be obtained by building a spanning tree of a given network \cite{scaman2017optimal}. For $m$-node centralized distributed architecture, the number of gradient oracle calls per each node is defined by  \eqref{SS} with $\sigma^2 = 0$, $L$ and $\mu$ correspond\dm{ing} to $\tilde{f}$.  The number of communication steps will be $d$ times \dm{more}, where $d$ is the distance between the origin (root) and farthest \dm{node}.} 
If we have only $q \ll m$ nodes, then we divide the data $\{f(x,\xi^k)\}_{k=1}^m$ \dm{into} $q$ blocks \dm{with} $l = m/q$ terms in each block. If $l$ is too large by itself on\dv{e} can reformulate \eqref{sum_rand} as follows \cite{mcmahan2016communication}
\begin{equation}\label{sum_rand_stoch}
\min_{x\in Q\subseteq \R^n}\tilde{f}(x):=\frac{1}{q}\sum_{k=1}^q \E[f_k(x,\eta^k)],   
\end{equation}
where $f_k(x,\eta^k) = f(x,\xi^{kl + \eta^k})$ and $\eta^k = i$ with probability $1/l$, $i=1,...,l$. Representation \eqref{sum_rand_stoch} allows to use bound \eqref{SS} in the stochastic case in a parallel manner at each node. 
The number of oracle calls per node also corresponds \g{(in general)} to \eqref{SS} and the number of communication steps is also $d$ times \dm{more} than \eqref{SS} with $\sigma^2=0$\g{.} 

Unfortunately, centralized architecture has a synchronization drawback and a high requirement for the master node \cite{scaman2017optimal}. To address these disadvantages to some extent, \dm{a} decentralized distributed architecture \dm{should be used} \cite{bertsekas1989parallel,kibardin1979decomposition}. \dm{This architecture} relies on two basic principles \sasha{\cite{nedic2020distributed}}:  every node communicate\dd{s} only with its neighbors, and  \ga{all communications are performed simultaneously}. The main difference here is \dm{a} simple strategy of communications: each node communicates only with all available direct neighbors. This architecture is more robust. In particular, it can be applied to time-varying (wireless) communication networks \cite{rogozin2018optimal}.

\dm{As}  problems \eqref{sum_rand} and \eqref{sum_rand_stoch} have \dm{ a definite structure of the sum type}, they can be solved much faster on one machine. For instance, using some incremental algorithms \cite{allen2016katyusha,lan2017optimal,lin2015universal,woodworth2016tight}, one can solve \eqref{sum_rand} $\sqrt{m}$ times cheaper in terms of the number of oracle calls, but not in terms of the number of iterations ( =communication steps). Unfortunately, this result prohibits  \g{parallelization}. Note, that for this problem in asynchronized mode ( at each step only two randomly chosen nodes can communicate )\dm{,} one can obtain such ($\sim\sqrt{m}$) an acceleration for the star-type communication network \cite{lan2018random}. Moreover, a `dual' analogue of this acceleration has  recently been proposed for \eqref{sum_rand}  \cite{hendrikx2018accelerated} and \eqref{sum_rand_stoch}  \ga{\cite{hendrikx2019accelerated,hendrikx2019asynchronous}} with arbitrary communication networks.

\dv{Before stating the contribution we  introduce the notions of condition number  ${\chi}$   for \gav{Laplacian} communication matrix of some  network  and the height of spanning tree denoted by $d$. Note, that $\sqrt{\chi}\ge d$ and typically $\sqrt{\chi} \le nd$ (see \cite{nedic2018graph}). The last bound corresponds to \dm{a} star topology \cite{gasnikov2017modern}  (the most simple centralized type architecture). In many interesting cases $\sqrt{\chi} = \tilde{O}(d)$ (see \cite{nedic2018graph,scaman2017optimal}).}

\subsection{Contribution}
\begin{itemize}
    \item {We justify the transition  from \dm{the} optimal centralized distributed complexity bounds for problems \eqref{sum_rand} and \eqref{sum_rand_stoch} in the smooth case to  decentralized ones\footnote{
 In the deterministic case this was partially done in \cite{li2018sharp}.}
by replacing $d$ 
 with $\sqrt{\chi}$, the average $L$ with the worse one and variance of  $f$ with the variance of $f_k$, that can be $m$ times \dm{more}.}\footnote{\ga{For instance, this takes place in the case when we have independent noise at each $f_k$. Note also that in this case in decentralized distributed optimization one can improve the variance dependence   and \dv{eliminate} factor $m$  (see \cite{olshevsky2019asymptotic,olshevsky2019non}). But, this is possible due to the worse estimate for the number of communication steps.}} \dv{Here and everywhere below  we keep  $\mu$ at the average level \dm{(without loss of generality, we can assume that each $f_k$ has the same $\mu$)} by using a trick from \cite{scaman2017optimal}.} \dv{ The announced results are also not improvable in terms of communication steps (rounds) \cite{arjevani2015communication,scaman2017optimal}.} 

By using different smoothing techniques \cite{allen2016optimal,nesterov2005smooth,scaman2018optimal}, we may \dm{lead} the non-smooth case to the smooth one with $L\sim 1/\e$. This allows to reduce the complexity estimate \eqref{SNS} by using \eqref{SS}. However, in general, this reduction \dm{makes the cost of  oracle calls more expensive}. Thus, we can only improve the communication steps (rounds) bound that corresponds to \eqref{SS} (up to a $\sqrt{\chi}$ factor) with $L\sim 1/\e$ and $\sigma^2 = 0$.\footnote{Note that such  tricks sometimes allow to obtain \dm{the} optimal (in terms of dependence on $\e$) communication round estimates \cite{arjevani2015communication,scaman2018optimal}.} Can we preserve the bound \eqref{SNS} for standard conception of oracle calls \dm{(primal oracle that gives $\nabla f_k$)} per node in decentralized approach by improving the number of communication steps? The answer is positive up to the replacement of \dd{the} average $M$ to \dd{the} wor\dd{st} one   \cite{lan2017communication,scaman2018optimal}. In our paper\dm{,} we simplify \dd{the} approaches \dd{proposed in} these articles to prove this result.
\end{itemize}
As a motivation for the second main result\dm{,}  we consider the problem of type 
\begin{equation}
\label{Sum}
\min_{x\in Q}f(x) := \frac{1}{m}\sum_{k=1}^m f_k(x).  
\end{equation}
where $f_k(x)$ has Fenchel--Legendre representation $f_k(x) = \max\limits_{y}\left\{\langle x,y \rangle - \vp_k(y)\right\}$ with convex $\vp_k(y)$. Such type of \avg{optimization problems arise, \dm{for instance}, in different applications of inverse problems, especially for linear problems in Hilbert space where we have to use discretization and sum-type functions naturally arise \cite{byrne2014iterative,gao2017distributed,gasnikov2017convex,vogel2002computational,ye2019optimization}.} \avg{The case of dual-friendly $f_k$ also arises} in the problem of Wasserstein barycenter calculation \cite{dvinskikh2019dual,dvurechenskii2018decentralize,uribe2018distributed}. Suppose that $\nabla \vp_k(y)$ is available but  $f_k(x), \nabla f_k(x)$ are unavailable. Moreover, sometimes  $\vp_k(y) = \E[\vp_k(y,\xi)]$ and it is worth to use $\nabla \vp_k(y,\xi)$ instead  of $\nabla \vp_k(y)$  \cite{dvinskikh2019dual,dvurechenskii2018decentralize}. \dv{Considering this example as one of the motivation for using dual oracle instead of primal one, we provide the second main result.}
\begin{itemize}
\item 
{We develop  optimal decentralized distributed algorithms with dual (stochastic) oracle for strongly convex objective in \eqref{Sum}.}
The approach is based on dual reformulation of \eqref{Sum} \cite{scaman2017optimal}. \dm{An} optimal algorithm for non-strongly convex  dual function with stochastic oracle was recently proposed in \cite{dvinskikh2019dual}. To propose an optimal method with stochastic dual oracle \dm{for} strongly convex \dm{primal objective,} \gav{we use} recent work \cite{foster2019complexity}. We notice  a rather unexpected result: we cannot improve (up to a logarithmic factor) the bound for the \dm{number of} dual stochastic gradient calculations in comparison with non-strongly convex dual objective.

We also notice that initially we \dd{were} motivated by \dm{the study} of the dual oracle not only as an application from \cite{dvinskikh2019dual,dvurechenskii2018decentralize,uribe2018distributed}. We also tried to find a simple explanation \gav{for} \dm {the} optimal communication step bounds \cite{arjevani2015communication,scaman2018optimal} in non-smooth case. One of the way\dm{s} to do it is Nesterov's dual smoothing technique \cite{nesterov2005smooth} that builds a bridge to the notion of dual oracle. This plan was partially (in the deterministic case) implemented in \cite{scaman2017optimal,uribe2018dual,uribe2018dual}. Here we generalize the results of these works for the stochastic dual oracle. 
\end{itemize}

\subsection{Paper organization}
The paper  is organi\dm{z}ed as follows. In Section \ref{sec1}\dm{,}  we propose optimal stochastic (parallelized) accelerated gradient methods for stochastic convex optimization problems. In Sections \ref{prim} and \ref{dual}\dm{,}  we apply the results of  Section \ref{sec1} to stochastic convex optimization problems with affine type of constraints (of type $Ax=0$). We describe the modern stochastic (parallelized) accelerated gradient methods which are optimal both in terms of (stochastic) oracle calls and matrix-vector multiplications $Ax$. In Sections \ref{prim}\dm{,}  we \dd{are} focusing on primal methods, in Section \ref{dual}\dm{,}  we present dual ones.  Section \ref{dec}\dm{,}  describes the distributed primal and dual formulation of \dm{the} finite-sum minimization problem, and  \dm{presents distributed algorithms}. In Section \ref{main}\dm{,}  we incorporate \dm{the} proposed distributed decentralized method to get \dm{the} optimal bounds for \dm{the} finite-sum minimization problem using primal or dual oracle. Finally, we discuss future work and possible extensions. We notice that  all proposed methods  are  optimal in terms of communication steps  and \g{in many cases in terms of (parallel stochastic) primal/dual oracle calls.} 


\section{Stochastic convex optimization}\label{sec1}
First-order methods for the optimization problem of minimizing a convex function $f$ on a simple convex set $Q$, e.g., 
    \begin{align}\label{eq:gen_prob}
        \min_{x\in Q \subseteq \R^n} f(x),
    \end{align}
play a fundamental role in modern problems arising in machine learning and statistics. The complexity of these methods is measured by the number of iterations or (and) the number of oracle calls. For a deterministic oracle\dm{,} this concept can be identified.
By the first-order {oracle}, we  mean  a {black-box} model that for a given input $x \in Q$, returns the vector $\nabla f(x)$.

\dm{We say that a function $f$  is  $M$-Lipschitz continuous  if \footnote{Here and below in such type of assumptions (especially in the case when $Q$ is unbounded) instead of $\forall x\in Q$ we may write $\forall x \in Q: \|x - x^*\|_2 \le 2R$ \cite{gasnikov2017modern} (analogously for $y$).}
 \[ \forall x \in Q \qquad  \|\nabla f(x)\|_2 \leq M.\]}
We say that function $f$  is  $L$-smooth or has $L$-Lipschitz continuous gradient if
\begin{equation*}
   \forall x,y\in Q \qquad \|\nabla f(y) - \nabla f(x)\|_2 \leq L \|y-x\|_2.
\end{equation*}
We also say that function $f$  is $\mu$-strongly convex if
\begin{equation*}
   \forall x,y\in Q \qquad f(y)\geq f(x) +\la\nabla f(x), y-x\ra + \frac{\mu}{2} \|y-x\|_2^2.
\end{equation*}

\begin{algorithm}[h]
\caption{Similar Triangles Method {\tt STM($L$,$\mu$,$x^0$)},  $Q = \mathbb{R}^n$}
\label{Alg:STM}   
 \begin{algorithmic}[1]
\Require $\tilde{x}^0=z^0=x^0$, number of iterations $N$, $\alpha_0 = A_0=0$, $L$, $\mu$ 
\For{$k=0,\dots, N$}
  \State Set $\alpha_{k+1} = \frac{1+A_{k}\mu}{2L} + \sqrt{\frac{1+A_{k}\mu}{4L^2}+\frac{A_k\left(1+A_k\mu\right)}{L}}$, $A_{k+1} = A_k + \alpha_{k+1}$
\State $\tilde{x}^{k+1} = (A_kx^k+\alpha_{k+1}z^k)/A_{k+1}$
\State $z^{k+1} = z^k - \frac{\alpha_{k+1}}{1+\dmdv{A_{k+1}}\mu} \left( \nabla f(\tilde{x}^{k+1}) \dmdv{+} \mu\dmdv{\left(z^k - \tilde{x}^{k+1}\right)}\right)$
\State $x^{k+1}=(A_kx^k+\alpha_{k+1}z^{k+1})/A_{k+1}$
  
\EndFor
        \Ensure    $ x^N$ 
\end{algorithmic}
 \end{algorithm}

Accelerated gradient methods (e.g., Algorithm~\ref{Alg:STM} ({\tt STM})
 \cite{gasnikov2018universal, nesterov2018lectures, Lan2019lectures}) allow to obtain the optimal number of iterations and number of gradient oracle calls for problem~\eqref{eq:gen_prob}  as described  in Table~\ref{T:deter}, where $R = \|x^0 - x^*\|_2$ is the   Euclidean distance from the starting point $x^0$ to the solution $x^*$ \dm{of \eqref{eq:gen_prob} that corresponds to the minimum of this norm}, and $\e$ is the desired precision in function value. 
\begin{remark}
\dv{{
For a composite optimization problem with composite term $h(x)$\dm{,}  step 4 of Algorithm \ref{Alg:STM} is replaced by the more general operator \cite{gasnikov2018universal, nesterov2018lectures}
 \begin{align*}
     z^{k+1} = \arg\min_{z\in Q} \left\{ \sum_{l=0}^{k+1}\alpha_l\left(
     \langle \nabla f(\tilde{x}^l),z - \tilde{x}^l \rangle +  h(z) + \frac{\mu}{2}\|z-\tilde{x}^l\|_2^2\right) \right.\\
     \left.+\frac{1}{2}\|z-\tilde{x}^0\|_2^2\right\}. 
 \end{align*}
  If $h(x)$ has $L_h$-Lipschitz gradient in the $\ell_2$-norm then due to Theorem 9 \cite{stonyakin2019gradient} and Theorem 19 \cite{stonyakin2019inexact}  
it suffices to solve auxiliary problem with accuracy (in terms of function value)
 \[O\left(\frac{(\alpha_{k+1}\e)^2(A_{k+1}\mu + 1)}{(A_{k+1}L_h R)^2}\right)\ge O\left(\frac{\e^3}{LL_h^2R^4}\right),\]
 where $\e$ is desired accuracy (in function value) for initial problem \eqref{eq:gen_prob}.\\
 If $\mu = 0$ we can also generalize this step for the non-Euclidean case and using restarts \cite{gasnikov2018universal}  generalize such a method on $\mu > 0$. Note, that by using restarts with {\tt STM($L$,0,$x^0$)} one can eliminate the gap  from $\ln(LR^2/\e)$ to $\ln(\mu R^2/\e)$ between lower bounds and the bounds for {\tt STM($L$,$\mu$,$x^0$)} without restarts \cite{gasnikov2018universal}. The same remains true for the stochastic \dm{oracle}.}}
 
\end{remark}

    	\begin{table}[h]
    	\caption {
    \dm{The} optimal number of first-order oracle calls (number of iterations $N$)
    	}
\label{T:deter}
\begin{center}
{\small
\begin{tabular}{ |c| c| c| c| c|}
 \hline
&{ \makecell{ $\mu${-strongly convex}\\ and $L$-smooth} } & { \makecell{ $L$-smooth}} &   { \makecell{ $\mu${-strongly convex}}} &  \\
 \hline
\makecell{\#iterations}& $\sqrt{\frac{L}{\mu}}\ln \left(\frac{\mu R^2}{\e}\right)$  &  $\sqrt{\frac{LR^2}{\e}}$ & $\frac{M^2}{\mu\e}$ & $\frac{M^2R^2}{\e^2}$    \\
 \hline
\makecell{\#oracle calls \\
of  $\nabla f(x)$
 }& $\sqrt{\frac{L}{\mu}}\ln \left(\frac{\mu R^2}{\e}\right)$  &  $\sqrt{\frac{LR^2}{\e}}$ & $\frac{M^2}{\mu\e}$ & $\frac{M^2R^2}{\e^2}$    \\
 \hline
\end{tabular}
}
\end{center}
\vskip -0.1in
\end{table}

Generally, iteration complexity \dm{is determined by the complexity of calculating the gradient, which  can be computationally expensive.} Thus, stochastic approximations of the true gradient can be used instead. In this case, or when the true gradient is unavailable (if e.g., function $f$ is given in the form of expectation $f(x):=\E[f(x,\xi)$]) we denote the {inexact} (or noise-corrupted) first-order {oracle} as $\nabla f(x, \xi)$, given by a blackbox model with stochasticity (noise) $\xi$ corrupting the true gradient. Assume that
\begin{equation}
\label{bias}
\left\|\E [\nabla f(x, \xi)] -\nabla f(x)\right\|_2\le\delta=O(\e/R)
\end{equation}
and
\begin{equation*}
\E\left[\exp\left({ \|\nabla f(x, \xi)-\E[\nabla f(x, \xi)]\|_2^2}{\sigma^{-2}}\right)\right]\le \exp(1),
\end{equation*}
then with probability at least $1 - \beta$, we have
$f(x^N) - f(x^*)\le \e$
after  
$$N=\min\left\{O\left(\sqrt{\frac{LR^2}{\e}}\right), 
 O\left(\sqrt{\frac{L}{\mu}}\ln\left(\frac{LR^2}{\e}\right)\right)\right\}$$
 iterations of { \tt STM} using the approximated gradient (instead of the real one $\nabla f(\tilde{x}^{k+1})$)
 \begin{equation}
 \label{batch_grad}
\nabla^{r_{k+1}} f(\tilde{x}^{k+1},\{\xi_i^{k+1}\}_{i=1}^{r_{k+1}}) = \frac{1}{r_{k+1}}\sum_{i=1}^{r_{k+1}} \nabla f(\tilde{x}^{k+1},\xi_i^{k+1}),
 \end{equation}
 where $\xi_1^{k+1}, \dots, \xi_{r_{k+1}}^{k+1}$ are i.i.d from the same distribution as $\xi$ and the batch size is
 \begin{equation}
 \label{batch_size}
 r_{k+1} = O\left(\frac{\sigma^2\alpha_{k+1}\ln(N/\beta)}{(1+A_{k+1}\mu)\e}\right).
 \end{equation}
Moreover, the total number of oracle calls\footnote{
Oracle calls can be easily and fully  parallelized (on $r_{k}$ processors) at each iteration. Note, that for $\nabla^{r_k} f(x,\{\xi_i\}_{i=1}^{r_k})$ we can reduce the variance $\sigma^2:=O(\sigma^2/r_k)$. 
} is (this bound is optimal up to logarithmic factors)
\begin{align*}
    \sum_{k=0}^N r_k = O(N) +\min&\left\{O\left(\frac{\sigma^2R^2}{\e^2}{\ln\left(\frac{\sqrt{LR^2/\e}}{\beta}\right)}\right),\right.\\
    &\left.O\left(\frac{\sigma^2}{\mu\e}{\ln\left(\frac{LR^2}{\e}\right)\ln\left(\frac{\sqrt{L/\mu}}{\beta}\right)}\right)\right\}.
\end{align*}
We refer to such a variant of {\tt STM} as {\tt BSTM($L$,$\mu$,$\sigma^2$,$x^0$)} (batched {\tt STM($L$,$\mu$,$x^0$)}).

Thus, using minibatches for constructing an approximation of the true gradient allows us to keep the optimal number of iterations for stochastic methods, as presented in Table~\ref{T:deter}, where we skip high probability logarithmic multipliers. The number of stochastic oracle calls for this case is shown in Table \ref{T:stoch}. 

\begin{remark}
\dm{We notice that in assumption \eqref{bias}, $R=\|x^0-x^*\|_2$. Generally, in such type of assumptions, $R$ is the diameter of $Q$ (see \cite{cohen2018acceleration,aspremont2008smooth}) (it is not a compact set when $Q = \mathbb{R}^n$).} To obtain such a generalization we have to use the advanced recurrent technique to bound $\|z^k - x^*\|_2$ from \cite{dvinskikh2019dual,gorbunov2018accelerated} and  \cite[Chapter 2]{gasnikov2017modern}. Further we provide the sketches how to get this result (for simplicity $\sigma = 0, \mu = 0$).\\
1. For \dm{an} inexact gradient $\tilde{\nabla} f(x)$ satisfying  for all $x,y$
\begin{align}
\label{model}
&f(x)+\langle\tilde{\nabla} f(x),y-x\rangle - \delta_1\avg{\|y-x\|_2} \le f(y) \nonumber \\ 
&\hspace{3cm}\le f(x)+\langle\tilde{\nabla} f(x),y-x\rangle + \frac{L}{2}\|y-x\|_2^2 +\delta_2,
\end{align}
{\tt STM} outputs\footnote{\sasha{Note, that according to \cite{poljak1981iterative} even for the last point of gradient descent on a simple quadratic optimization problem we cannot guarantee convergence without proper stopping rule. With proper stopping rule in \eqref{bias} it is required (see  \cite[Theorem 7, item 6.1.3]{polyak1987introduction}) $\delta \sim \e^2$ that is worse than what we have $\delta \sim \e$. But we can guarantee standard convergence of noisy gradient descent under \eqref{bias} in (Cesaro) average \cite{gasnikov2017modern} (not for the last point). The results below generalize \cite{gasnikov2017modern} on a proper accelerated method (STM).}} $x^N$ such that \sasha{\cite{devolder2014first,dvinskikh2020accelerated}}
\[f(x^N) - f(x^*) = O\left(\frac{LR^2}{N^2} + \delta_1\avg{\tilde{R}} + N\delta_2\right),\]
\avg{where $\max\left\{\|\tilde{x}^k - x^*\|_2,\|z^k - x^*\|_2,\|x^k - x^*\|_2\right\}\le \tilde{R}$.}\\
2. Since
\[\langle\tilde{\nabla} f(x) - \nabla f(x),y-x\rangle \le \frac{1}{2L}\|\tilde{\nabla} f(x) - \nabla f(x)\|_2^2 + \frac{L}{2}\|y-x\|_2^2,\]
one can consider $\delta_2 := \delta^2/(2L)$ and $L:=2L$ in \eqref{model}.\\
\avg{3.} In the deterministic case\dm{,}    $\tilde{R} = R$  \sasha{with the proper stopping rule of the algorithm (see  \cite[formulas (2.17), (2.18) in Chapter 2]{gasnikov2017modern} and \cite{gasnikov2018universal})}. In the stochastic case \sasha{(with high probability)} $\tilde{R} = O(R)$ (see \cite{dvinskikh2019dual,gorbunov2018accelerated}) for {\tt STM} with the proper batch-size~\eqref{batch_size}.
\end{remark}

In particular, for the case of non-smooth objective, the stochastic oracle does not yield gains compared to its deterministic counterpart.

	\begin{table}[H]
\caption { 
\dm{The} optimal number of {stochastic} (unbiased) first-order oracle calls 
}
\label{T:stoch}
\begin{center}
{\small
\begin{tabular}{ |c |c |c |c |c| }
 \hline
& \makecell{ $\mu$-strongly convex \\ and $L$-smooth  }  &   \makecell{  $L$-smooth } &   \makecell{ $\mu$-strongly convex } &  \\
 \hline
\makecell{\#iterations}& $\sqrt{\frac{L}{\mu}}\ln \left(\frac{\mu R^2}{\e}\right)$  &  $\sqrt{\frac{LR^2}{\e}}$ & $\frac{M^2}{\mu\e}$ & $\frac{M^2R^2}{\e^2}$    \\
 \hline
\makecell{\#oracle calls\\ of $\nabla f(x,\xi)$ } &\makecell{$\max\left\{  \frac{\sigma^2}{\mu\e},\right.$\\
$\left.  \sqrt{\frac{L}{\mu}}\ln \left(\frac{\mu R^2}{\e}\right) \right\}$}  &  \makecell{$\max\left\{\frac{\sigma^2R^2}{\e^2},\sqrt{\frac{LR^2}{\e}}\right\}$} & $\frac{M^2}{\mu\e}$ & $\frac{M^2R^2}{\e^2}$   \\
 \hline
\end{tabular}
}
\end{center}
\vskip -0.1in
\end{table}

\dm{In Table~\ref{T:stoch}, in the appropriate cells,  we assumed that the following inequalities hold:} $\mathbb{E}\|\nabla f(x, \xi) - \nabla f(x) \|_2^2 \leq \sigma^2$ and   $\E\|\nabla f(x,\xi)\|^2_2 \leq M^2$.\\

\dv{Both in  Table  \ref{T:deter} and \dm{in} Table~\ref{T:stoch}, the last two columns can be obtained from the corresponding first columns by choosing $L = M^2/(2\delta)$, where $\delta = \e/N$ (see \cite{gasnikov2018universal}). This is the idea of universal accelerated methods \cite{nesterov2015universal}, but with predefined $L$. Here and in all further tables we skip numerical constants.}

\section{Primal methods for stochastic convex optimization with affine constraints}\label{prim}
To build the complete theory of distributed primal and dual method we need to generalize the result of Tables~\ref{T:deter} and \ref{T:stoch} for the convex optimization problem\footnote{In decentralized optimization $A$ is taken to be $\sqrt{W}$ (square root of the Laplacian matrix of the communication network).}
\begin{equation}
\label{PP2}
\min_{\substack{Ax=0, \\ x\in Q}} f(x),    
\end{equation}
where $A\dm{\succeq} 0$ and $\text{Ker} A \neq \varnothing$. 
The purpose of this section is to develop such algorithms for \eqref{PP2} that are optimal in terms of the number of $\nabla f(x)$ calculations and the number of $A^T Ax$ calculations. In this section we use  Euclidean proximal setup \dm{\cite{ben-tal2001lectures}}. This is the only section where we  significantly rely on Euclidean prox-structure.

Denote by $R_y = \|y^*\|_2$  the $\ell_2$-norm of the smallest solution $y^*$ of dual (up to a sign) problem \eqref{DP}. Solution $y^*$ is not unique since $\text{Ker} A \neq \varnothing$.
From \cite{lan2017communication} we have such a bound
\begin{equation}
\label{R}
R_y^2 \le \frac{\|\nabla f(x^*)\|_2^2}{\lambda_{\min}^{+}(A^T A)}.
\end{equation}

Using the penalty method we rewrite \eqref{PP2} as follows
\begin{equation}
\label{penalty}
F(x) = f(x) + \frac{R_y^2}{\e}\| Ax\|_2^2 \to \min_{x\in Q}.  
\end{equation}
Next we use  \cite[Remark 4.2]{gasnikov2017modern} and get if the following holds
\[F(x^N)-\min_{x\in Q}F(x)\le \e,\]
then 
\[f(x^N)-\min_{x\in Q\ga{, Ax = 0 }}f(x)\le \e, \quad \|Ax^N||_2\le \ga{\frac{(1+\sqrt{5})\e}{2R_y}}.\]

We start with the smooth case and assume that $Q = \mathbb{R}^n$. If $f$ has $L$-Lipschitz continuous gradient 
then we can solve \gav{\eqref{penalty}} by {\tt STM} (or {\tt BSTM} in the stochastic case) considering the second term to be composite \cite{gasnikov2018universal,nesterov2013gradient}. In this case, we obtain \dm{the} optimal number of $\nabla f(x)$ (or $\nabla f(x,\xi)$) calculations, see Tables~\ref{T:deter},~\ref{T:stoch}. But the total number of $A^T Ax$ calculations is
 \[\tilde{O}\left(\sqrt{\lambda_{\max}(A^T A)/\lambda_{\min}^{+}(A^T A)}\right)\] times more
\dv{{ since  $\text{Im} A = \text{Im} A^T = (\text{Ker} A)^{\perp}$ and $Q = \mathbb{R}^n$ and as a consequence of these facts the auxiliary problem can be divided into two subproblems: minimization  of quadratic form with matrix of the form $(R_y^2/\e)\gav{A^T A} + cI$  ($c$ is some positive constant and $I$ is identity matrix) on $(\text{Ker} A)^{\perp}$ and minimization  of quadratic form with matrix of the form $cI$ on $\text{Ker} A$. Linear terms do not play any role in complexity. The complexity of the auxiliary problem \dm{is} determine\dm{d} by the worst (reduced on corresponding subspace) conditional number of these two subproblems. Obviously, the first one is worse. The reduced the conditional number is \[\frac{\lambda_{\max}\left((R_y^2/\e)\gav{A^T A} + cI\right)}{\lambda_{\min}^{+}\left((R_y^2/\e)\gav{A^T A} + cI\right)} \le \frac{\lambda_{\max}(\gav{A^T A})}{\lambda_{\min}^{+}(\gav{A^T A})}.\]}}
This factor arises because of the complexity of the auxiliary problem. We refer to these approaches as {\tt PSTM} and {\tt PBSTM} (Penalty {\tt STM} and {\tt BSTM}). Here and below we skip arguments of the algorithms if they are obvious from the context.

In non-smooth case ($f$ is $M$-Lipschitz)\dm{,} 
we use the Sliding algorithm \cite{lan2016gradient}, \ga{\cite{Lan2019lectures}} . If $\mu = 0$ according to \cite{lan2016gradient} this algorithm requires (see Tables~\ref{T:deter},~\ref{T:stoch} for comparison)
\begin{center}
$O\left(\sqrt{\frac{\lambda_{\max}(A^T A)R_y^2 R_x^2}{\e^2}}\right)$ 
\end{center}
calculations of $A^T Ax$ and
\begin{center}
$O\left(\frac{M^2R_x^{2}}{\e^2}\right)$ 
\end{center}
calculations of $\nabla f(x)$, where $R_x = \|x^0 - x^*\|_2$.

If  we have unbiased $\nabla f(x,\xi)$ with $\sigma^2$-sub-Gaussian variance \cite{jin2019short} instead of $\nabla f(x)$, i.e.,
\begin{equation*}
\E\left[\exp\left({ \|\nabla f(x, \xi) - \nabla f(x)\|_2^2}{\sigma^{-2}}\right)\right]\le \exp(1),
\end{equation*}
with $\sigma^2 = O(M^2)$ (for compact notation\footnote{\dv{In general $M^2$ is replaced by $ M^2 + \sigma^2$ in the stochastic case.}}), then
\ga{the bound for calculations of $A^T A$ does not change and the bound for calculations of $\nabla f(x,\xi)$ is the same as it was for the number of calculations of $\nabla f(x)$ in deterministic case}
(up to a logarithmic high-probability deviations factor).

By using a restart technique  \cite{uribe2018dual} we can generalize this method for $\mu$-strongly convex $f$: 
\[O\left(\sqrt{\frac{\lambda_{\max}(A^T A)R_y^2 }{\mu\e}}\ln\left(\frac{\mu R_x^2}{\e}\right)\right)\]
calculations of $A^T Ax$ and
\[O\left(\frac{M^2}{\mu\e}\right)\]
calculations of $\nabla f(x)$ ($\nabla f(x,\xi)$).

We  call this approach by {\tt R}-Sliding (Restart Sliding).

\section{Dual methods for stochastic convex optimization with affine constraints}\label{dual}

Now we assume that we can build a dual problem for
\begin{equation}
\label{PP}
\min_{\substack{Ax=0,\\ x\in Q}} f(x),    
\end{equation}
where $\text{Ker} A \neq \varnothing$.

\begin{remark}
\dv{We notice that turning to a dual problem does not oblige us using dual oracle. Instead\dm{,} we can use a primal oracle and the Moreau theorem \cite{Rockafellar2015} with Fenchel-Legendre representation. This maximization problem can be solved  using the first-order oracle for the function $f$. \gav{But such an approach does not allow to obtain \dm{the} optimal bounds on the  number of primal first-order oracle calls.}
Note  that typically in decentralized optimization $A$ in \eqref{PP} is taken \dm{as the} square root of \dm{the} Laplacian matrix $W$ of the communication network \cite{scaman2017optimal}. But  in the  asynchronized case the square root $\sqrt{W}$ replaced by incidence matrix $M$ \cite{hendrikx2018accelerated} ($W = M^T M$). Then in asynchronized case instead of accelerate\dm{d} methods for \eqref{DP} one should use an accelerated (block) coordinate descent method \cite{dvurechensky2017randomized,gasnikov2017modern,hendrikx2018accelerated,shalev-shwartz2014accelerated}.}
\end{remark}

\begin{algorithm}[h]
\caption{{\tt PDSTM}}
\label{Alg:PDSTM}   
 \begin{algorithmic}[1]
\Require $\tilde{y}^0=z^0=y^0=0$, number of iterations $N$, $\alpha_0 = A_0=0$
\For{$k=0,\dots, N$}
  \State Set $\alpha_{k+1} = \frac{1}{2L_{\dm{\psi}}} + \sqrt{\frac{1}{4L_{\dm{\psi}}^2}+\frac{A_k}{L_{\dm{\psi}}}}$, $A_{k+1} = A_k + \alpha_{k+1}$
\State $\tilde{y}^{k+1} = (A_ky^k+\alpha_{k+1}z^k)/A_{k+1}$
\State $z^{k+1} = z^k - \alpha_{k+1} \nabla \psi(\tilde{y}^{k+1}) = z^k - \alpha_{k+1} A x(A^T\tilde{y}^{k+1})$
\State $y^{k+1}=(A_ky^k+\alpha_{k+1}z^{k+1})/A_{k+1}$

\EndFor
        \Ensure    $y^N$, $x^N = \frac{1}{A_N}\sum_{k=0}^N \alpha_k x(A^T\tilde{y}^k)$
\end{algorithmic}
 \end{algorithm}

The dual problem (up to a sign) is \dm{following}
\begin{align}\label{DP}
\psi(y) = \vp(A^T y) &= \max_{x\in Q}\left\{\langle y,Ax\rangle - f(x)\right\} = \langle y,Ax(A^T y)\rangle - f(x(A^T y)) \nonumber\\
&=  \langle A^Ty,x(A^T y)\rangle - f(x(A^T y)) \to\min_{y}.    
\end{align}
If $f$ is $\mu$-strongly convex in the $\ell_2$-norm, then $\psi$ has   $L_{\psi}=\frac{\lambda_{\max}(A^T A)}{\mu}$--Lipschitz continuous gradient in the $\ell_2$-norm\footnote{Here and below we can also consider \gav{other} norms (see \cite{uribe2018dual} for details).} \cite{kakade2009duality,Rockafellar2015}. In this case we can apply {\tt STM($L{\gav{_\psi}}$,0,0)} to \eqref{DP}. Note that due to Demyanov--Danskin's theorem $\nabla \psi(y) = Ax(A^T y)$ \cite{Rockafellar2015}. Similarly to \cite{anikin2017dual,chernov2016fast} one can prove that
\begin{align}
\label{PD}
f(x^N)-f(x^*)&=f(x^N)-f(x(A^T y^*))\le f(x^N) + \psi(y^N) \nonumber \\
& = O\left(\frac{L_{\psi}R_y^2}{N^2} \right), \quad \|Ax^N\|_2 = O\left(\frac{L_{\psi}R_y}{N^2} \right),
\end{align}
where $R_y = \|y^*\|_2$ is the radius of solution of \eqref{DP} which is  the smallest in the $\ell_2$-norm, see \eqref{R}.  
We  call this approach by {\tt PDSTM} (Primal-Dual {\tt STM}).

If we  have only a stochastic (randomized) unbiased  model $\nabla \vp(\lambda, \xi)|_{\lambda = A^T y} = x(A^T y,\xi)$ with $\sigma^2_{\vp}$\gav{-sub-Gaussian} variance, i.e. 
\begin{align*}
    &\E\left[\exp\left({ \| \nabla \vp(\lambda, \xi) - \nabla \vp(\lambda)\|_2^2}{\sigma^{-2}_{\vp}}\right)\right] \\
    &\hspace{3cm}= \E\left[\exp\left({ \| x(A^T y, \xi)- x(A^T y)\|_2^2}{\sigma^{-2}_{\vp}}\right)\right]
    \le \exp(1),
\end{align*}
then for {\tt BSTM($L_{\psi}$,0,$\sigma^2_{\psi}$,0)} where $\sigma^2_{\psi} = \lambda_{\max}(A^TA)\sigma^2_{\vp}$ with probability $\ge 1 - \beta$~\eqref{PD} holds true \cite{dvinskikh2019dual}. We refer to this algorithm as {\tt SPDSTM} (Stochastic {\tt PDSTM}).\\

In the case when $\psi$ \dm{from} \eqref{DP} is additionally $\mu_{\psi}$-strongly convex in the $\ell_2$-norm in\footnote{Since $\text{Im} A  = (\text{Ker} A^T)^{\perp}$ we will have that all the points $\tilde{y}^k,z^k,y^k,$ generated by STM and \dm{methods based on STM}, belong to $y^0 + (\text{Ker} A^T)^{\perp}$. That is, from the point of view of estimates this means, that we can consider $\psi$ to be $\mu_{\psi}$-strongly convex everywhere.} $y^0 + (\text{Ker} A^T)^{\perp}$ (if $f$ has $L$-Lipschitz gradient in the $\ell_2$-norm and $Q = \mathbb{R}^n$ then $\mu_{\psi} = \lambda_{\min}^{+}(A^T A)/L$ \cite{kakade2009duality,Rockafellar2015}, where $\lambda_{\min}^{+}(A^T A)$ is the minimal positive eigenvalue of $A^T A$) we need to use another approach. Because of primal-duality \dm{\cite{nesterov2009primal, nemirovski2010accuracy}}  we have to put $\mu_{\psi} = 0$ in {\tt STM} and \dm{in methods based on STM}  \gav{({\tt STM($L_{\psi}$,$\mu_{\psi}$,$y^0$)} is not primal-dual method when  $\mu_{\psi} > 0$)}. The restart technique (see, e.g. \cite{gasnikov2017modern}) also does not work here because in~\eqref{PD} we have to use  $R_y = \|y^0\|_2 + \|y^0 - y^*\|_2$ in general. That is why we take here $y^0 = 0$. So the main trick here is the following relation \cite{allen2018make,anikin2017dual,nesterov2012make}
\begin{equation}
\label{grad_norm}
f(x(A^T y))- f(x^*) 
\le\langle\nabla \psi(y), y\rangle = \langle Ax(A^T y), y\rangle.    
\end{equation}
From \eqref{grad_norm}, \dm{to satisfy}
\begin{equation*}
\ga{f(x^N)-f(x^*)=f(x(A^T y^N))-f(x(A^T y^*))
\le \ga{2}\e,
\|Ax^N\|_2 \le \e/R_y,}
\end{equation*}
 it is sufficient to find such $y^N$ \ga{($\|y^N\|_2\le 2R_y$)} that 
\begin{equation}
\label{GN}
\|\nabla \psi(y^N)\|_2\le \e/R_y.
\end{equation}

Recently, there appear accelerated methods with the proper rate of convergence in terms of the norm of the gradient {\tt OGM-G} \cite{gasnikov2017modern,kim2018optimizing}:
\[\|\nabla \psi(y^N)\|_2 = O\left(\frac{L_{\dm{\psi}}\ga{\|y^0 - y^*\|_2}}{N^2}\right)\dm{ = O\left( \frac{L_{\psi} \|\nabla \psi(y^0)\|_2}{\mu N^2} \right)}.\]
After $\bar N =  O\left(\sqrt{\frac{L_{\psi}}{\mu_{\psi}}}\right)$ iterations of { \tt OGM-G} we  have 
\[\|\nabla \psi(y^{\bar{N}})\|_2 \le \frac{1}{2}\|\nabla \psi(y^0)\|_2.\]
So after  $l = \log_2\left(\gav{\|}\nabla \psi(y^{0})\|_2 \frac{R_y}{\e}\right)$ restarts ($y^0 : = y^{\bar{N}}$) we have
\eqref{GN}. We  denote such an approach  by {\tt ROGM-G} (Restart {\tt OGM-G}). This approach requires $$O\left(\sqrt{\frac{L_{\psi}}{\mu_{\psi}}}\ln\left(\|\nabla \psi(y^{0})\|_2 \frac{R_y}{\e}\right)\right)$$  
of $\nabla \psi(y)$ (that is  $Ax(A^T y)$) calculations. \dv{The key inequality  to prove this fact is 
$$\|y^0 - y^*\|_2^2 \le \frac{1}{\mu^2_{\psi}}\|\nabla \psi(y^0)\|_2^2.$$} \dm{This holds due to  $$\frac{\mu}{2}\|y^0-y^*\|^2_2 \leq \psi(y^0) - \psi(y^*) \leq \frac{1}{2\mu_{\psi}} \|\nabla \psi(y^0)  \|_2^2.$$}

The same result with the replacement \[\dm{\sqrt{\frac{L_{\psi}}{\mu_{\psi}}}}\ln\left(\|\nabla \psi(y^{0})\|_2 \frac{R_y}{\e}\right) \to \dm{\sqrt{\frac{L_{\psi}}{\mu_{\psi}}}}\ln\left(2L_{\psi}^2\frac{R^4_y}{\e^2}\right)\] can be obtained by using {\tt STM($L_{\psi}$,$\mu_{\psi}$,0)} with \dm{ bound $\dm{\sqrt{\frac{L_{\psi}}{\mu_{\psi}}}}\ln \left(\frac{L_\psi R_y^2}{\e'} \right)$ (see \cite{nesterov2010introduction}}) and desired accuracy $\e' = \frac{\e^2}{2L_{\psi}R^2_y}$. This follows from 
$$\frac{1}{2L_{\psi}}\|\nabla \psi(y^N)\|_2^2 \le \psi(y^N) - \psi(y^*) \dm{\leq \e'}.$$

Now  \dm{we} consider {\tt RRMA+AC-SA$^{2}$}  \cite{foster2019complexity} (see also \cite{allen2018make} in the non-accelerate, but composite case). This algorithm converges as follows (for simplicity we skip polylogarithmic factors and high probability terminology) 
$$\|\nabla \psi(y^N)\|_2^2 = \tilde{O}\left(\frac{L^2_{\psi}\|y^0 - y^*\|_2^2}{N^4} + \frac{\sigma^2_{\psi}}{N} \right) = \tilde{O}\left(\frac{L^2_{\psi}\|\nabla \psi(y^0)\|_2^2}{\mu^2_{\psi} N^4} + \frac{\sigma^2_{\psi}}{N} \right).$$
\dm{We assume that  at each iteration,  $\nabla \psi(y,\xi)$ with sub-Gaussian variance $\sigma^2_{\psi}$ is available} \cite{jin2019short} (see also above). If we use restarts with the size of each restart $\bar{N} = \tilde{O}\left(\sqrt{\frac{L_{\psi}}{\mu_{\psi}}}\right)$ (see above) and use batched gradient \eqref{batch_grad} with batch size (at $k$-th restart; $\bar{y}^{k}$ is the output point from the previous restart)
$$r_{k+1} = \tilde{O}\left(\frac{\sigma^2_{\psi}}{\bar{N}\|\nabla \psi(\bar{y}^{k+1})\|_2^2}\right).$$ 
then $\|\nabla \psi(\bar{y}^l)\|_2 \le \e/R_y$ after $l = O\left(\log_2 \left(\|\nabla \psi(y^0)\|_2R_y/\e\right)\right)$ restarts. Therefore, the total number of oracle calls is
$$\tilde{O}\left(\frac{\sigma^2_{\psi}R^2_y}{\e^2}\right).$$ 
Note that the same bound 
takes place in the non-strongly convex case ($\mu_{\psi} = 0$). From \cite{allen2018make,jin2019short} it is known that this bound cannot be improved. However, we may expect that this bound can be reduced to  $\tilde{O}(\sigma^2_{\psi}/(\mu_{\psi}\e))$ ( see Table~\ref{T:stoch} for stochastic primal oracle).  For the stochastic dual oracle\dm{,} such a reduction is probably  impossible.
We call this approach by {\tt R-RRMA+AC-SA$^{2}$} (Restart {\tt RRMA+AC-SA$^{2}$}).

\section{Decentralized distributed optimization}\label{dec}

Now we show how to \dm{present} \eqref{Sum} in a decentralized distributed manner
\begin{equation}\tag{P1}\label{P1}
  \min_{x\in Q\subseteq \R^n} f(x) := \frac{1}{m}\sum_{k=1}^m f_k(x).   
\end{equation}

This particular representation of the objective in \eqref{P1} allows involving distributed methods which are particularly necessary for large-scale problems handling the large quantities of data and which are based on the idea of agents' cooperative solution of the global problem \cite{bertsekas1989parallel}. For a given  multi-agent network system\dm{,} we privately assign each function $f_k$ to the agent $k$  and suppose that agents can exchange the information with their neighbors (e.g., send and receive vectors). We define this system through the Laplacian matrix $\bar W \in \R^{m\times m}$ of some graph (communication network) $G=(V,E)$ with the set $V$ of $m$ vertices and the set of edges $E = \{(i,j): i,j \in V\}$ as follows
{\small	\begin{align*}
	 \bar{W}_{ij} = \begin{cases}
	-1,  & \text{if } (i,j) \in E,\\
	{\rm deg}(i), &\text{if } i= j, \\
	0,  & \text{otherwise,}
	\end{cases}
	\end{align*}}
 where ${\rm deg}(i)$ is the degree of vertex $i$ (i.e., the number of neighboring nodes).
	 
 From the definition of the matrix $\bar W$ it can be easily seen that $\bar W$ establishes the communication of agents and allows only the communication between neighboring nodes. Moreover, due to connectivity of graph $G$  the vector $\bld{1}_m = (1,...,1)^T\in\R^m$ is the unique (up to a scaling factor) eigenvector of $\bar W$ associated with the eigenvalue $\lm =0$, which allows us to compactly rewrite  the consensus agreement $x_1= ... = x_m\in \R^n$  as  $W\x=0$, moreover, as $\sqrt{W}\gav{\x}=0$ (see \cite{scaman2017optimal}), where $W \gav{= \bar{W}} \otimes I_n$ is the Kronecker product of the Laplacian matrix $\bar W\in \R^{m}$ and the identity \sasha{(unit)} matrix $I_n$ and $\x = [x_1^T, ..., x_m^T]^T \in \R^{mn}$.

To present  problem \eqref{P1} in a distributed fashion we rewrite it with introducing the artificial consensus equality constraints and then change these constraints to  one affine constraint with the  communication matrix $W$ as follow\dm{s}

\begin{equation*}
\min_{\substack{x_1= \dots =x_m , \\ x_1,\dots, x_m \in Q \subseteq \R^n }} F(\x) := \frac{1}{m}\sum_{k=1}^m f_k(x_k)
\end{equation*}
or in another form
\begin{equation}\tag{P2}\label{P2}
  \min_{\substack{\sqrt{W}\x=0, \\ x_1,\dots, x_m \in Q \subseteq \R^n }} F(\x) := \frac{1}{m}\sum_{k=1}^m f_k(x_k),
\end{equation}
where 
all $f_k$ are $M$-Lipschitz, $L$-smooth and $\mu$-strongly convex (it is possible that, $L = \infty$
or (and) $\mu = 0$).

We also consider the stochastic version of  problem (P2), where $f_k(x_k) = \E[f_k(x_k,\xi_k)]$. 
We consider the unbiased stochastic primal oracle \dm{that} returns  $\nabla f_k(x_k, \xi_k)$ (where  $\xi = \{\xi_k\}_{k=1}^m$ are independent) under the following $\sigma^2$-sub-Gaussian variance condition (for all $k=1,...,m$)
\begin{equation*}
\E\left[\exp\left(\frac{ \|\nabla f_k(x_k, \xi_k)- \nabla f_k(x_k)\|_2^2}{\sigma^2}\right)\right]\le \exp(1).
\end{equation*}

Problem \eqref{P2} can be considered to be a particular case of problem \eqref{PP2} with the following replacements
 \begin{algorithmic}[1]
\Statex \textbullet \hspace{0.3cm} $A = \sqrt{W}$
\Statex \textbullet \hspace{0.3cm} $L_F = {L}/{m}$
\Statex \textbullet \hspace{0.3cm} $\mu_F = {\mu}/{m}$
\Statex \textbullet \hspace{0.3cm} $\|\nabla F(\x)\|_2^2 \le M_F^2 = {M^2}/{m}$
\Statex \textbullet \hspace{0.3cm} $\sigma_F^2 = O\left({\sigma^2}/{m}\right)$
\Statex \textbullet \hspace{0.3cm} $R_{\x}^2 = \|\x^0 - \x^*\|_2^2 = m\|x^0 - x^*\|_2^2 = mR^2$
\Statex \textbullet \hspace{0.3cm} $R_{\y}^2 = \|\y^*\|_2^2 \le {\|\nabla F(\x^*)\|_2^2}/{\lambda_{\min}^{+}(W)}\le{M^2}/{(m\lambda_{\min}^{+}(W))}$
\end{algorithmic}

The main observation in the primal approach (see Section~\ref{prim}) is as follows \cite{scaman2017optimal}:
\begin{center}
\fbox{{
$A^T Ax  = Wx$  (calculated in a decentralized distributed manner)}
}
\end{center}

If  each function $f_k$ is a dual-friendly  \dm{(dual function is available calculated by the Fenchel--Legendre transform \cite{uribe2018dual})}
then we can construct the dual problem to problem \eqref{P2} with dual Lagrangian variables   $\y = [y_1^T \in \R^n,\cdots,y_m^T \in \R^n]^T \in \R^{mn}$  

\begin{equation}\label{D2}\tag{D2}
\min_{\y \in \R^{mn}}\Psi(\y) := \ga{\frac{1}{m}}\Phi(\ga{m}\sqrt{W}\y):= \frac{1}{m}\sum_{k=1}^{m}
\vp_k(\ga{m}[\sqrt{W}\y]_k), 
\end{equation}
where $\vp_k(\lambda_k) = \max\limits_{x_k\in Q \subseteq \R^n} \{\la \lambda_k,x_k\ra  - f_k(x_k)\}$ and the
vector $[\sqrt{W}\x]_k$ represents the $k$-th $n$-dimensional block of $\sqrt{W}\x$.

We also consider the stochastic version of  problem \eqref{D2}, where $\vp_k(\ga{\lambda}_k) = \E[\vp_k(\ga{\lambda}_k,\xi_k)]$. We consider the unbiased stochastic dual oracle returns  $\nabla \vp_k(\ga{\lambda}_k, \xi_k)$ (where  $\xi = \{\xi_k\}_{k=1}^m$ are independent) under the following $\sigma_{\vp}^2$-sub-Gaussian variance  condition (for all $k=1,...,m$)
\begin{equation*}
\E\left[\exp\left({ \|\nabla \vp_k(\ga{\lambda}_k, \xi_k)- \nabla \vp_k(\ga{\lambda}_k)]\|_2^2}{\sigma_{\vp}^{-2}}\right)\right]\le \exp(1).
\end{equation*}

Problem \eqref{D2} can be considered as a particular case of problem \eqref{DP} with 
 \begin{algorithmic}[1]
\Statex \textbullet \hspace{0.3cm} $A = \sqrt{W}$
\Statex \textbullet \hspace{0.3cm} $\sigma_{\ga{\Psi}}^2 = O\left(\ga{\lambda_{\max}(W)m\sigma_{\vp}^2}\right)$
\end{algorithmic}

The main observation in the dual approach (see Section~\ref{dual}) is as follows \cite{scaman2017optimal}:
since $x(A^T y) = \x(\sqrt{W}\y)$ we should change the variables as follows
 \begin{algorithmic}[1]
\State   \hspace{0.3cm} $\tilde{\y} :=\sqrt{W}\tilde{\y}$
\State \hspace{0.3cm} $\mathbf{z} := \sqrt{W}\mathbf{z}$ 
\State  \hspace{0.3cm} $\y := \sqrt{W}\y$
\end{algorithmic}

It is obvious that {input}, {output} and steps 3--5 of Algorithm~\ref{Alg:PDSTM} are changed such that they can be \dm{performed} in a decentralized distributed manner. For that we  just multiply the corresponding steps \dm{by} $\sqrt{W}$.

 \section{Main Results}\label{main}
 In this section, we present the rates of convergence for problems \eqref{P1} and \eqref{D2} (and their stochastic counterparts) in terms of the number of iterations (communication steps) and the number of (parallelized) oracle calls.
For the primal problem\dm{,} we present the results to achieve $\e$-precision in objective residuals, and for the dual problem we seek to achieve $\e$-precision in  duality gap or  \dm{primal} objective residuals (in smooth strongly convex case). Feasibility constrains are smaller than  $\e/R_\y$. 

For brevity, we introduce the condition number of the Laplacian matrix $W$ as  follows
\begin{equation}\label{chi}
\chi = \frac{\lambda_{\max}(W)}{\lambda_{\min}^{+}(W)}, 
\end{equation}
where $\lambda_{\min}^{+}(W)$ is the minimal positive eigenvalue of $W$, and $\lambda_{\max}(W)$ is the maximal eigenvalue of $W$.
Now we are ready to present our main results incorporated in multiple tables: Tables \ref{tab:DetPrimeOr}-- \ref{T:dist_stoch}. Th\dm{ese} results are obtained by direct substitution of constants from Section \ref{dec} to the problems from Sections~\ref{prim} and \ref{dual}.

    	\begin{table}[ht!]
\caption { The optimal bounds for primal deterministic oracle}
\label{tab:DetPrimeOr}
\begin{center}
{\small
\begin{tabular}{|c| c| c| c| c|}
 \hline
 & \makecell{ $f_k$ is $\mu$-strongly\\  convex\\ and $L$-smooth} 
 &   \makecell{$f_k$ is $L$-smooth} 
 & \makecell{ $f_k$ is $\mu$-strongly\\ convex} &  \\
 \hline
 \makecell{\#communication \\ rounds} 
 & \makecell{$\widetilde O\left(\sqrt{\frac{L}{\mu}\chi} \right)$}   
 & \makecell{ $ \widetilde O\left({\sqrt{\frac{LR^2}{\e} \chi}}\right)$ }
 &  \makecell{  $O\left(\sqrt{\frac{M^2}{\mu\e}\chi} \right)$} 
 &  \makecell{  $O\left( \sqrt{\frac{M^2R^2}{\e^2}\chi} \right) $}   \\
 \hline
\makecell{\#oracle calls of \\ $\nabla f_k(x_k)$\\ per node $k$} 
& $ \widetilde O\left(\sqrt{\frac{L}{\mu}} \right)$ 
&  $O\left(\sqrt{\frac{LR^2}{\e}}  \right) $ 
&   $O\left(\frac{M^2}{\mu\e} \right) $ 
&   $O\left( \frac{M^2R^2}{\e^2}\right)$\\
 \hline
Algorithm & \makecell{{\tt PSTM}, \\ $Q=\R^n$} & \makecell{{\tt PSTM},\\ $Q=\R^n$} &  {\tt R}-Sliding &  Sliding\\
 \hline
\end{tabular}
}
\end{center}
\end{table}

\begin{table}[ht!]
\caption {The optimal bounds for primal stochastic (unbiased) oracle} 
\label{T:stoch_prima_oracle}
{\small
\hspace{-1cm}\begin{tabular}{ |c| c| c| c| c|}
 \hline
&{\makecell{ $f_k$ is 
$\mu$-strongly \\ convex\\ and $L$-smooth} } &  {\makecell{ $f_k$ is  $L$-smooth} }&  {\makecell{ $f_k$ is 
$\mu$-strongly \\ convex} } &  \\
 \hline
\makecell{\#communication \\
rounds}& \makecell{ $\widetilde O\left(\sqrt{\frac{L}{\mu}\chi} \right)$}   & \makecell{ $ \widetilde O\left({\sqrt{\frac{LR^2}{\e} \chi}}\right)$ }&  \makecell{  $O\left(\sqrt{\frac{M^2}{\mu\e}\chi} \right)$} &  \makecell{  $O\left( \sqrt{\frac{M^2R^2}{\e^2}\chi} \right) $}  \\
 \hline
   \makecell{\#oracle calls \\ of $\nabla f_k(x_k, \xi_k)$\\
 per node $k$} & {\makecell{ $ \widetilde O\left(\max\left\{\frac{\ga{\sigma}^2}{\mu\e},
  \sqrt{\frac{L}{\mu}}\right\}\right)$ }}  & {\makecell{$ O\left(\max\left\{\frac{\ga{\sigma}^2R^2}{\e^2}, \sqrt{\frac{LR^2}{\e}}  \right\}\right)$}} & { {\makecell{  $O\left(\frac{M^2\ga{+\sigma^2}}{\mu\e} \right)$}}} &{\makecell{  $O\left(\frac{(M^2\ga{+\sigma^2})R^2}{\e^2}\right)$}} \\
 \hline
Algorithm & \makecell{{\tt PBSTM},\\ $Q=\R^n$}& \makecell{{\tt PBSTM},\\ $Q=\R^n$} & \makecell{ Stochastic \\{\tt R}-Sliding} & \makecell{ Stochastic \\ Sliding} \\
 \hline
\end{tabular}
}
\end{table}

    	\begin{table}[ht!]
\caption {The optimal bounds for dual deterministic oracle}
\label{tab:distrDetCOm}
\begin{center}
{\small
\begin{tabular}{ |c |c |c| }
 \hline
 & \makecell{ $f_k$ is $\mu${-strongly convex}\\ and $L$-smooth} 
&  \makecell{ $f_k$ is $\mu$-strongly convex }  \\
 \hline
\makecell{\#communication  \\ rounds} & \makecell{ $\widetilde O\left(\sqrt{\frac{L}{\mu}\chi} \right)$ } &  \makecell{ $O\left(\sqrt{\frac{M^2}{\mu\e}\chi}\right) $ }  \\
 \hline
\makecell{  \#oracle calls of\\ $\nabla \vp_k(\lambda_k)$
 per node $k$} & \makecell{ $\widetilde O\left(\sqrt{\frac{L}{\mu}\chi} \right)$ } &  \makecell{ $O\left(\sqrt{\frac{M^2}{\mu\e}\chi}\right) $ }  \\
 \hline
Algorithm & \makecell{{\tt ROGM-G} or {\tt STM}, $Q = \R^n$}  &
 \makecell{ {\tt OGM-G} or {\tt PDSTM} }  \\
 \hline
\end{tabular}
}
\end{center}
\end{table}

	\begin{table}[ht!]
\caption{The optimal bounds for dual stochastic (unbiased) oracle}
\label{T:dist_stoch}
\begin{center}
{\small
\begin{tabular}{ |c| c| c| }
 \hline
&{\makecell{ $f_k$ is $\mu${-strongly convex}\\ and $L$-smooth} }  & \makecell{ $f_k$ is $\mu$-strongly convex } \\
 \hline
\makecell{\#communication \\ rounds} & \makecell{ $\widetilde O\left(\sqrt{\frac{L}{\mu}\chi} \right)$ }  & \makecell{ $O\left(\sqrt{\frac{\ga{M^2}}{\mu\e}\chi}\right) $ } \\
 \hline
\makecell{\#oracle calls of \\
$\nabla \vp_k(\lambda_k,\xi_k)$
 per node $k$}  &  \makecell{$\widetilde O\left(\max\left\{\pd{\frac{M^2\sigma_{\ga{\vp}}^2}{\e^2}\chi}, \sqrt{\frac{L}{\mu}\chi} \right\}\right)$ }  & \makecell{ $O\left(\max\left\{ \frac{M^2\sigma_{\vp}^2}{\e^2}\chi, \sqrt{\frac{ M^2}{\mu\e}\chi} \right\}\right)$ }   \\
 \hline
Algorithm & {\tt R-RRMA+AC-SA$^{2}$}, $Q=\R^n$ &  {\tt SPDSTM}  \\
 \hline
\end{tabular}
}
\end{center}
\end{table}

Note that the bounds on communication steps (rounds) are optimal (up to a logarithmic factor) due to \cite{arjevani2015communication,scaman2017optimal,scaman2018optimal}.
Bounds for the oracle calls per node are probably optimal \ga{in the class of  methods with  optimal number of communication steps} (up to a logarithmic factor) in the deterministic case \ga{\cite{allen2018make,foster2019complexity,woodworth2018graph}} and optimal \g{for the non-smooth stochastic primal oracle and stochastic dual oracle}
for parallel architecture.\footnote{In parallel architecture the bounds on stochastic oracle calls per node of type $\max\{B,D\}$ can be \g{parallel} up to $B/D$ processors.} For stochastic oracle the bounds hold  in terms of high probability deviations (we skip the corresponding logarithmic factor).

We  emphasize that the difference between centralized (or parallel) estimates and \g{obtained} decentralized ones is not only in the replacement of $d$ by $\tilde{O}(\sqrt{\chi})$ in the smooth cases for the primal oracle \avg{and the meaning of $L$}. In the stochastic smooth strongly convex case (one can also consider the convex case)   we know that the total number of primal oracle calls \cite{kulunchakov2019estimate1,kulunchakov2019estimate2,kulunchakov2019generic,lan2018random} is \[\tilde{O}\left(m + \sqrt{m\frac{L}{\mu}} + \frac{\sigma^2}{\mu\e}\right).\] This bound is optimal but it uses an incremental oracle and does not  imply  full parallelization. The best known way to parallelize it is described in \cite{lan2018random}. For full parallelization one should use a standard accelerated scheme without variance reduction and incremental oracle \cite{woodworth2018graph}. In this case, another bound for the total number of oracle calls occurs, that is
\[\tilde{O}\left(m \sqrt{\frac{L}{\mu}} + \frac{\sigma^2}{\mu\e}\right).\]
But this bound assumes the natural way of parallelization or centralized distribution of calculations. In the last case for a graph of diameter $d$ with  $m$  nodes we  have the following number of oracle calls per node
\[\tilde{O}\left( \sqrt{\frac{L}{\mu}} + \frac{\sigma^2}{m\mu\e}\right) = \tilde{O}\left(\max\left\{\frac{\sigma^2}{m\mu\e}, \sqrt{\frac{L}{\mu}}\right\}\right)\]
 and the following number of communication steps 
 \[\tilde{O}\left(d \sqrt{\frac{L}{\mu}} \right).\]
 For decentralized architecture (see Table~\ref{T:stoch_prima_oracle}) the number of oracle calls per node and the number of communication steps are
 \begin{equation}\label{OptimalBound}
\tilde{O}\left(\max\left\{\frac{\sigma^2}{\mu\e}, \sqrt{\frac{L}{\mu}}\right\}\right) \text{ and } \tilde{O}\left( \sqrt{\frac{L}{\mu}\chi} \right).   
 \end{equation}
 respectively. Unfortunately, the factor $m$  is no longer presented in $\sigma^2$ in the decentralized case. It is interesting to note, that it is possible to propose such a decentralized distributed algorithm that requires
 \[O\left( \frac{\sigma^2}{m\mu\e}\right)\] oracle calls per node (stochastic gradients calculations) \cite{olshevsky2019asymptotic,olshevsky2019non}. However, this algorithm is not optimal in terms of communication steps. \avg{Moreover, to the best of our knowledge, it is an open question whether is \eqref{OptimalBound} optimal bound in terms of oracle calls (per node) in the class of methods with the optimal number of communication steps.}

Note \ga{also} that the \pd{blue} bound in Table~\ref{T:dist_stoch} seems to be rather unexpected at first sight for us. \ga{But we 
\sasha{hypothesize} that this bound is optimal not only in terms \dm{of the number} of communication steps but also  in terms of \dm{the number of} oracle calls (per node) in the class of  methods with the optimal number of communication steps.}

\ga{The detailed proofs of the statements collected in this paper takes more than 90 pages. These can be found in the arXiv preprint \cite{gorbunov2019optimal}:}
\begin{center}
   \avg{\url{https://arxiv.org/pdf/1911.07363.pdf}.}
\end{center}

\section{Discussion}
Below we \dm{outline various areas for} further work.
\begin{itemize}
    \item  We can expect that the  results can be improved by
    replacing
    \dd{the} first-order methods with primal and dual deterministic oracle by tensor methods ($p=2,3$) from \cite{nesterov2018implementable}. \ga{However, for the moment we do not know any such results.} \ga{F}or the dual approach we \ga{also do not} know how to use the trick $A=\sqrt{W}$ (see \cite{scaman2017optimal}). Here we should take $A = W$, which (with additional increased complexity of auxiliary problem) make\dd{s} the bounds on communication steps wor\gav{se}.  The basic fact in the dual approach is the following. To solve auxiliary problems we have to calculate the values of the form $\nabla^p_y \vp(Wy)$ on different vectors \cite{carmon2016gradient,nesterov2018implementable,nesterov2018lectures}.  \ga{T}his can be done multiplying $W$ by vectors (communications) and multiplying corresponding (block) diagonal tensor ( $\nabla^p_{\lambda} \vp(\lambda)|_{\lambda=Wy}$) by vectors (can be distributed among nodes)
    \item  The primal approach in the smooth case can be generalized for the (stochastic inexact) gradient-free oracle. The number of communication steps remain\dd{s} the same. The number of oracle calls become\dd{s} $\sim n$ times larger \cite{gorbunov2018accelerated,dvurechensky2017randomized}. \g{In the non-smooth case gradient-free (stochastic) decentralized distributed algorithm was developed in \cite{beznosikov2019derivative}}
    \item In \ga{\cite{hendrikx2019accelerated},} \cite{hendrikx2019asynchronous} \dm{asynchronized} distributed optimization was considered via dual accelerated (block) coordinate descent algorithms. The  primal approach proposed above allows asynchronized generalizations in the smooth case. For that we should use the (block) coordinate  version of {\tt STM}  \cite{dvurechensky2017randomized} and additional randomization of sum type when  $[Wx]_i$ is calculated. This will increase the number of communication steps $\sim\sqrt{n} \div n$ times
    \item Most of the results of this paper can be generalized to composite problems \cite{nesterov2013gradient}. \dd{Perhaps,} it is possible to make \dd{the} next step and try to generalize these results to more general types of models \cite{stonyakin2019gradient,stonyakin2019inexact}
    \item  For smooth convex centralized distributed optimization problems there exists a universal way to accelerate non-accelerated (stochastic, asynchronized etc.) algorithm  Catalyst \cite{lin2015universal}. The basic idea is using  a non-accelerated centralized distributed algorithm for the inner problem arising at each step of the Catalyst procedure
    \item  
    \dd{Perhaps,} it is 
    possible to generalize the primal approach described above on time-varying graphs \avg{\cite{rogozin2019projected}}. Moreover, these generalizations can be done also for the smooth stochastic  case 
    \item \ga{It seems the result of \cite{rogozin2018optimal} can be improved by using mixed communication: many decentralized steps alternate with centralized ones. In this case, one can use non-accelerated distributed decentralized algorithms, which are robust on time-varying graphs \cite{rogozin2018optimal}, and then accelerate them by using the Catalyst technique \cite{lin2015universal} and centralization. Since the graph is changing we should recalculate spanning tree \dm{whenever} we apply a centralized step. We expect that this mixed communication will be useful also for tensor schemes in decentralized optimization}
    \item \ga{The main scheme in the primal approach is based on the result formulated directly after \eqref{penalty}. This result \dd{does not} depend on convexity of the objective. So it \dd{would} be interesting to apply this scheme for non-convex distributed optimization problems \cite{sun2018distributed}}
\end{itemize}

\avg{Since the first version of this paper was submitted on arXiv, there appeared alternative explanations (for smooth problems) of the results for the primal deterministic oracle \sasha{\cite{fallah2019robust,kovalev2020optimal,li2020revisiting,xu2019accelerated,hendrikx2020dual}} and the primal stochastic oracle (strongly convex case) \cite{fallah2019robust}. Moreover, in \g{\cite{rogozin2020towards,ye2020multi}} (see also \cite{rogozin2019projected} for the non-accelerated case and \cite{scaman2018optimal} for lower bounds) for  the primal deterministic oracle (strongly convex case) it was shown that $L: = \max_{k} L_k$ used in this paper can be improved to $L: = L_f$, where $L_f$ is the Lipschitz gradient constant of \eqref{Sum}  (which can be much smaller \cite{tang2019practicality}). The same holds true for $\mu$. More interestingly, we expect that combinations of \g{\cite{rogozin2019projected,rogozin2020towards,ye2020multi}} allows to develop accelerated decentralized distributed primal algorithms on time-varying graphs (accelerated in $L$ and $\e$, but not in $\chi$). Moreover, based on  \g{\cite{rogozin2019projected,rogozin2020towards,ye2020multi}} we may expect that the answer to the open problem posed in the Section~\ref{main} is negative. \dm{The bound can be decreased by a factor $m$} for a stochastic primal oracle. The reasons for that are almost the same that we have for the $\max_{k} L_k \to L_f$ improvement.}

\avg{Based on the recent works \cite{koloskova2020unified,woodworth2020local,woodworth2020minibatch}  we describe below a hypothesis about \dm{the} optimal bounds in a more general situation. We consider the general sum-type problem (for simplicity we consider the case $Q = \mathbb{R}^n$, but it can be naturally generalized)
\begin{equation}
\label{GSP}
\min_{x\in \mathbb{R}^n} f(x) := \frac{1}{m}\sum_{k=1}^m f_k(x) = \frac{1}{m}\sum_{k=1}^m \E[f_k(x,\xi)].    
\end{equation}
We assume that all $f_k(x,\xi)$ in \eqref{GSP} satisfy
\begin{equation}\label{Lips}
\|\nabla f_k(y,\xi) - \nabla f_k(x,\xi)\|_{\sasha{2}} \le L \|y- x\|_2.
\end{equation}
Also we introduce\footnote{\avg{Note that nowadays it's quite popular to obtain estimates on rate of convergence that depend on the constants defined at the solution point $x^*$ \cite{stonyakin2020mirror}.}}
$$\sasha{\bar{\zeta}^2 } = \frac{1}{m}\sum_{k=1}^m \|\nabla f_k(x^*)\|_2^2$$
and
$$\sasha{\bar{\sigma}^2} = \frac{1}{m}\sum_{k=1}^m \E\left[\|\sasha{\nabla} f_k(x^*,\xi) - \nabla f_k(x^*)\|_2^2\right].$$
Assume now that at each iteration $t$ we may call one \sasha{time} an oracle (that returns \sasha{an independent realization of} $\nabla f_k(x^t,\xi^\sasha{t,k})$) at each node and make no more than one communication step with (in general, random) communication matrix $\bar{W}_k$. Moreover, we assume that \cite{koloskova2020unified}\footnote{\sasha{In \cite{koloskova2020unified} it was also assumed that symmetric matrix $\bar{W}$ determined by doubly stochastic matrix $P$: $\bar{W} = P - I$, where $I$ -- unit matrix. So do we. But below, when we generalize the results of \cite{koloskova2020unified} on a large class of algorithms that assumes more than one communication on one iteration, we may consider $\bar{W}$ to be the same as in section~\ref{dec} and $\chi$ is determined in \eqref{chi} style.}}
\sasha{
$$\E_{\bar{W}_l,...,\bar{W}_{l+\tau}}\left[\left\|(I+\bar{W}_{l+\tau})\cdot...\cdot(I+\bar{W}_l) x - \frac{1}{m}\bld{1}_m \bld{1}_m^T x \right\|_2^2\right] \le$$ $$\le (1 - \chi^{-1})\left\| x - \frac{1}{m}\bld{1}_m \bld{1}_m^T x\right\|_2^2.$$}
For the  class of algorithms described above  \sasha{the} required number of iterations $N$ to achieve an accuracy $\E[f(x^N)] - f(x^*)\le \varepsilon$ for the best known (for the moment) algorithms is}
\begin{equation*}\label{Koloskova1}
\tilde{O}\left(\left(\frac{ LR^2}{\varepsilon}\right)^{\sasha{\alpha}}+\frac{\bar{\sigma}^2R^2}{m\varepsilon^2} + \pd{\frac{\sqrt{ L}R^2\left(\bar{\zeta}\chi\tau + \bar{\sigma}\cdot(\chi\tau)^{\beta}\right)}{\varepsilon^{3/2}}(1-\chi^{-1})^{I[\tau = 1]}} + \left(\frac{\chi LR^2}{\varepsilon}\right)^{\sasha{\alpha}}\sasha{\tau}\right),
\end{equation*}
\avg{where $\beta \in [0,1]$, \sasha{$\alpha = 1 \text{ or } 1/2,$} \sasha{here} $I$  is a  function such that $I[\rm true] = 1$, $I[\rm false] = 0$, and when $f_k$ is $\mu$-strongly convex}
\begin{equation*}\label{Koloskova2}
\tilde{O}\left(\left(\frac{ L}{\mu}\right)^{\sasha{\alpha}}+ \frac{\bar{\sigma}^2}{m\mu\varepsilon} + \pd{\frac{\sqrt{ L}\left(\bar{\zeta}\chi\tau + \bar{\sigma}\cdot(\chi\tau)^{\beta}\right)}{\mu\sqrt{\varepsilon}}(1-\chi^{-1})^{I[\tau = 1]}} + \left(\frac{\chi L}{\mu}\right)^{\sasha{\alpha}}\sasha{\tau}\right).
\end{equation*}
\avg{This bound with $\beta = 1, \sasha{\alpha = 1}$ \sasha{(non-accelerated case)}  was obtained in \cite{koloskova2020unified} for simple  \sasha{(local)} decentralized SGD. Roughly speaking, $\beta = 0$ corresponds to the lower bound for this algorithm \sasha{\cite{koloskova2020unified,yuan2020federated,woodworth2020local,woodworth2020minibatch}}. \sasha{In \cite{karimireddy2019scaffold} by using a variance reduction trick in federated learning architecture ($\bar{W}_{p\tilde{\tau} + q} \equiv 0$ for $q=1,...,\tilde{\tau} - 1$ and $\bar{W}_{p\tilde{\tau}}$ is a full communication matrix, hence $\chi = 1$, $\tau = \tilde{\tau}$),  the SCAFFOLD algorithm was proposed with $\alpha = 1$ and (under some additional assumptions) without middle (\pd{blue}) terms. Also in the federated learning architecture in case $f_1 = ... = f_m$ it seems that middle (\pd{blue}) terms and the factor $\tau$  in the last terms of the described bounds can be eliminated under some natural additional assumptions \cite{woodworth2020local}. Here we assume at least one communication. If there is exactly one communication it will take place at the very end. This result means that in a federated learning setup under $f_1 = ... = f_m$ the frequency of communications $\tau$ does not play any significant role, which was previously mentioned in \cite{godichon2020rates}.}
\sasha{In} centralized architecture it seems that it is possible to obtain \sasha{additional} acceleration \sasha{of  the algorithms mentioned above ($\alpha = 1/2$, but the transition $\chi \to \sqrt{\chi}$ in the last terms requires additional assumptions)} by using proper accelerated envelopes like Catalyst \cite{Barre2020principle,dvinskikh2020accelerated,ivanova2019adaptive,kulunchakov2019generic} with the precision $\varepsilon'=\tilde{O}\left(\varepsilon\sqrt{\max\{\varepsilon/(LR^2),\mu/L}\}\right)$ (in the function) of the solution of auxiliary problem at each outer iteration. \sasha{It is highly likely that this is true for more general architectures, see examples in \cite{hendrikx2020dual,li2020revisiting}.}}

\avg{If we \sasha{remove the condition that at each iteration we can only call stochastic oracle one time and make no more than one communication,}
then the} \pd{blue} 
\avg{\sasha{term}  
can be eliminated \sasha{with $\alpha = 1/2$}. More precisely, the first two terms correspond to the total number of oracle calls per node and the last \sasha{term} 
-- to the number of communications steps.  Up to a factor $m$ in the denominator we have developed these results in the paper, but with \sasha{$\tau = 1$ and} fixed $W$.\footnote{\sasha{For time-varying communication graphs in the eneral case for the moment it seems that we should put the factor $\chi$ instead of $\sqrt{\chi}$ in the last term \cite{rogozin2020penalty_}.}}} \sasha{We  also note that it is an open problem  whether it is possible in the general (accelerated) situation \sasha{( which differs from the \cite{koloskova2020unified})} to determine $\bar{\zeta}^2$ and $\bar{\sigma}^2$ only at point $x^*$? \sasha{It seems, that for the current moment of time we have a positive answer only with respect to $\bar{\zeta}^2$.}}

\avg{The following generalizations are related with the case
$$f_k(x) = \frac{1}{l}\sum_{j=1}^{l} f_{kj} (x).$$
In this case with additional assumptions about proximal and dual friendly \sasha{$f_k$} it is possible to reduce worth case constant $L$ to the <<average>> one \cite{hendrikx2020optimal}. In \cite{hendrikx2020optimal} this looks like a variance reduction acceleration, but the nature of the effect (also explained in \cite{hendrikx2020optimal}) based on coordinate descent acceleration \cite{nesterov2017efficiency} for the dual problem formulation.  \sasha{In \cite{hendrikx2020dual}  was proposed a dual-free generalization of \cite{hendrikx2020optimal} with a bit worse oracle complexity estimate. The main idea is to apply a non-accelerated coordinate descent for the dual problem with Bregman divergence determined by the dual function itself.} \g{In \cite{li2020optimal} an optimal algorithm both for communications steps and oracle calls per node was developed.}
}

\avg{Another way to obtain \dm{a} better result for the required number of communication steps is possible if $\{f_{kj}\}$ have i.i.d. nature (that is typically for data science applications). In this case $f_k$ is statistically similar to $f$. Based on this fact in centralized architecture\dm{,} we may use on a master node statistical preconditioned algorithms that can significantly reduce the required number of communications with slaves \cite{hendrikx2020statistically}.}

\section*{Acknowledgment}
The work of D. Dvinskikh in \dm{Sections 1 and 5} \g{was supported by RFBR 19-31-51001 and in Section 6}  was funded by \dm{the} Russian Science Foundation (project 18-71-10108). 
 The work of A. Gasnikov was supported by the Ministry of Science and Higher Education of the Russian Federation (Goszadaniye) no. 075-00337-20-03, \avg{project No. 0714-2020-0005.} 

We would like to thank  F.~Bach, \avg{P.~Dvurechensky,} \ga{E.~Gorbunov, \sasha{A.~Koloskova,} A.~Kulunchakov, J.~Mairal,  A.~Nemirovski, A.~Olshevsky, S.~Parsegov,} B.~Polyak\ga{,} N.~Srebro \ga{ A.~Taylor and C.~Uribe} for useful \ga{discussions}.

\bibliographystyle{apalike}

\bibliography{all_refs3}

\end{document}